\documentclass[12pt,openany,brazil,reqno]{article}
\usepackage{amsmath,amsthm,amsfonts,amssymb,amscd}
\usepackage{graphics}
\usepackage[latin1]{inputenc}
\usepackage[portuges]{babel}

\renewcommand{\thefootnote}

\def\leaderfill{\leaders\hbox to .8em{\hss .\hss}\hfill} 
\def\_#1{{\lower 0.7ex\hbox{}}_{#1}}

\begin{document}

\centerline{\large\bf{Cohomology-free diffeomorphisms on tori}}

\vglue .4in

\centerline{\bf Nathan M. dos Santos}

\medskip

\centerline{\bf Instituto de Matemática}

\centerline{\bf Universidade Federal Fluminense}

\centerline{\bf 24020-605, \,\, Niterói - Rio de Janeiro - Brazil}

\vglue .5in

\noindent{\bf Abstract}.\, We study cohomology-free $(c.f.)$ diffeomorphisms of the torus $T^n$. A diffeomorphism is $(c.f)$ if every smooth function $f$ on $T^n$ is cohomologous to a constant $f_0$ i.e. there exists a $C^\infty$ function $h$ so that $h-h\circ \varphi = f-f_0$\,. We show that the only $(c.f.)$ diffeomorphisms of $T^n$ are smooth conjugations of Diophantine translations. This is part of a conjecture of A. Katok [H, Problem 17].

\bigskip

\noindent \textbf{1.\,Introduction.}\, A diffeomorphism $\varphi\colon T^n \to T^n$ is given on the covering $R^n$ by $\widetilde{\varphi} = A+F+\alpha$ where $A$ is an integer $n \times n$ matrix with $\det A = \pm 1$ and $F(x+p) = F(x)$ for all $p \in \mathbb{Z}^n$ and $\alpha \in R^n$, $\alpha \ne 0$.

A translation $T_\alpha\colon T^n \to T^n$ is {\it Diophantine\/} if
\begin{equation}
||k\cdot\alpha|| \ge \frac{C}{|k|^{n+p}}\quad C,\beta > 0 \tag{1.1}
\end{equation}
where $k = (k_1,\dots,k_n) \in \mathbb{Z}^n-\{0\}$, $k\cdot\alpha + k_1\,\alpha_1 +\cdots+ k_n\cdot\alpha_n$\,,
\begin{equation}
||x|| = \inf\{(x-\ell), \ell \in \mathbb{Z}^n\},\, |x| = {\operatornamewithlimits{\sup}}_j|x_j| \,\,\, [AS] \tag{1.2}
\end{equation}

A $(c.f.)$ diffeomorphism leaves invariant a volume form and it is uniquely ergodic and minimal.

\bigskip

\noindent\textbf{Proposition 1.}\, \textit{There exists an invariant volume form for a cohomology-free diffeomorphism.}

\medskip

\noindent\textbf{Proof:}\, Let $\Omega_0$ be a volume form on $M$. Thus
\begin{equation}
\varphi^*\Omega_0 = \det\, D\varphi\Omega_0 \tag{1.3}
\end{equation}
since $\varphi\colon T^n \to T^n$ is $(c.f.)$ then there exist a constant $c$ and a $C^\infty$ function $h\colon M \to R$ such that
\begin{equation}
\log|\det\,D\varphi| = h - h \circ \varphi + c. \tag{1.4}
\end{equation}
We will show that $c=0$. Consider the volume form
\begin{equation}
\Omega = \exp\, h\Omega_0\,. \tag{1.5}
\end{equation}
Now from (1.4) and (1.5) we have
\begin{align*}
|\varphi^*\Omega| = (\exp\,h\circ\varphi)|\varphi^*\Omega_0| &= (\exp\,h\circ\varphi)|\det\,D\varphi|\Omega_0\\
&= \exp(h\circ\varphi+\log|\det\,D\varphi|)\Omega_0\\
&= \exp(h+c)\Omega_0\\
&= \exp(c)\exp\,h\Omega_0 \tag{1.6}\\
&= \exp(c)\Omega\,\text{by}\, (1.5)
\end{align*}
thus
\begin{equation}
\left\vert\int_M \varphi^*\Omega\right\vert |\deg(\varphi)| = (\exp\,c) \int_M\Omega \tag{1.7}
\end{equation}
from (1.3) and (1.4) we may assume that $\displaystyle{\int}_M \Omega=1$ by choosing a convenient function $h\colon M \to R$.

Now from (1.7) we get
\begin{equation}
|\deg(\varphi)| = \exp(c) = 1 \tag{1.8}
\end{equation}
thus $c=0$ and from (1.5) and (1.6) we have
$$
\varphi^*\Omega = \deg(\varphi)\Omega.
$$
\rightline{$\square$}

\bigskip

\noindent\textbf{Proposition 2.}\, \textit{The entropy of a $(c.f.)$ diffeomorphism $\varphi\colon T^p \to T^p$ vanishes.}

\noindent\textbf{Proof:} By Proposition 1, 
\begin{equation}
\log|\det\,D\varphi| = h\circ\varphi-h. \tag{2.1}
\end{equation}
Thus
$$
\int_{T^p} \log|\det\,D\varphi|d\mu = 0
$$
since 
$$
\int_{T^p}(h\circ\varphi-h)d\mu = 0
$$
\pagebreak

and by Ruelle inequality
\begin{equation}
h_\mu(\varphi) \le \int_{T^p} \log|\det\,D\varphi|d\mu = 0 \tag{2.2}
\end{equation}
\rightline{[M2] \qquad\qquad $\square$}

\vglue .2in

\noindent\textbf{Proposition 3.}\, \textit{Any power of a $(c.f.)$ diffeomorphism $\varphi\colon T^p \to T^p$ is also $(c.f.)$.}

\medskip

\noindent\textbf{Proof:} Let $\varphi\colon T^p \to T^p$ be a $(c.f.)$ diffeomorphism then $\varphi^r\colon T^p \to T^p$, \, $\forall\,r \in \mathbb{Z}^+$ is $(c.f.)$. We have to show that the equation
\begin{equation}
h\circ \varphi^r-h = f, \quad \forall\, f \in C_0^\infty (T^p),\, r \in \mathbb{Z}^+ \tag{3.1}
\end{equation}
has a unique solution $h \in C_0^\infty(R^p) \Longleftrightarrow$ the sequence
\begin{equation}
S_n(\varphi^r)f = f + f \circ \varphi^r +\cdots+ f \circ (\varphi^r)^n, \, n \in \mathbb{Z}^+ \tag{3.2}
\end{equation}
is uniformly bounded i.e.
$$
||S_n(\varphi^r)|| < C, \,\, C > 0
$$
by [F, Lemma 52].
\rightline{$\square$}

\bigskip

\noindent\textbf{Theorem.}\, \textit{Let $\varphi\colon T^p \to T^p$ be a $(c.f.)$ diffeomorphism. Then $\varphi$ is conjugate to a Diophantine translation $\tau_\alpha\colon T^p \to T^p$.}

\medskip

\noindent\textbf{Proof:} For, by Proposition 2 above the entropy of $\varphi$ is zero and by [M1] the spectral radius $sp(\varphi_*)=1$ thus all eigenvalues of $\varphi_*$ are roots of unity and by the Lefschetz fixed point theorem 1 is an eigenvalue. By Proposition 3 there exists $r \in \mathbb{Z}^+$ such that 1 is the only eigenvalue of $\varphi_*^r$ and $\varphi^r$ is $(c.f.)$. By [SL2, Corollary 1.7] $\varphi^r$ is conjugate to the Diophantine translation $\tau_\alpha^r\colon T^p \to T^p$ by a diffeomorphism homotopic to the identity $\psi\colon T^p \to T^p$. Thus $\psi$ conjugate $\varphi$ to the Diophantine translation $\tau_\alpha$\,.

\rightline{$\square$}

\newpage

\centerline{\bf References}

\begin{itemize}
\item[{[AS]}] J.L. Arraut and N.M. dos Santos. Linear foliations on $T^n$.\, Bol. Soc. Brasil. Mat. 21 (1991), 189-204.
\item[{[SL1]}] R.U. Luz and N.M. dos Santos,\, Minimal homeomorphisms on Low-dimension tori.\, Ergod. Th. \& Dynam. Sys. (2009), 1515-1521.
\item[{[SL2]}] N.M. dos Santos, R. Urzúa Luz,\, Cohomology-free diffeomorphism of low-dimension tori.\, Ergod. Th. \& Dynam. Sys. (1998), 18, 985-1006.
\item[{[H]}] S. Hurden,\, Problems on rigidity of group actions and cocycles. \, Ergod. Th. \& Dynam. Sys. (1985), 65-76.
\item[{[F]}] H. Furstenberg, \, Strict ergodicity and transformation of the torus.\, Amer. J. of Math. 83 (1961), 573-601.
\item[{[M1]}] A. Manning, \, Topological entropy and the first homology group, \, Warwick Lecture.\, Notes in Mathematics 468 (1975).
\item[{[M2]}] Mañé,\, Ergodic Theory and Differentiable Dynamics.\, Springer-Verlag, 1987.
\end{itemize}

\end{document}